\documentclass[10pt]{article}

\usepackage{amssymb}
\usepackage{amsmath}
\usepackage{setspace}
\usepackage{eufrak}

\usepackage{theorem}

\newtheorem{lem}{Lemma }[section]
\newtheorem{defi}[lem]{Definition}
\newtheorem{theo}[lem]{Theorem}
\newtheorem{prop}[lem]{Proposition}

\newtheorem{fait}[lem]{Fact}
\newtheorem{cor}[lem]{Corollary}

\newcommand{\proof}{\noindent{\bf Proof}\par\nobreak}
\newcommand{\qed}{\ \hfill$\square$\smallskip}

\newcommand{\<}{\langle}

\renewcommand{\>}{\rangle}

\newcommand{\Rc}{{\rm{RC}}}

\title{Note on the Cantor-Bendixson rank of limit groups}
\author{Abderezak OULD HOUCINE}

\date{}

\begin{document}
\maketitle
\begin{abstract}
We show that the Cantor-Bendixson rank  of a limit group is finite as well as  that  of a limit group of a linear group.
\end{abstract}

\bigskip
Let $T$ be a
universal theory.  If $T$ has at most countably
many finitely generated models,  then one can associate  to each
finitely generated model $\cal M$ of $T$ an ordinal  rank, denoted
$Rk(\mathcal M)$, as defined in \cite{ould-f-g}. We reformulate
the definition of $Rk$ in the context of topological spaces. Let $X$
be a compact Hausdorff  topological space. Let $D(X)=\{x \in X | x
\hbox{ is not isolated in }X\}$. We define inductively
$D^\alpha(X)$ on  ordinals as follows:

\smallskip
$\bullet$ $D^0(X)=X$,

$\bullet$ $D^{\alpha+1}(X)=D(D^\alpha(X))$,

$\bullet$ $D^\alpha(X)= \bigcap_{\beta<\alpha}D^\beta(X)$, for
$\alpha$ a limit ordinal.

\bigskip
The following fact can be easily extracted from \cite{kechris}. 
\begin{fait}\label{fact1} There exists a least  ordinal $\alpha$ such that
$D^\alpha(X)=D^\beta(X)$ for any $\beta>\alpha$. If $X$ is
separable then $D^\alpha(X)=\emptyset$ or
$|D^\alpha(X)|=2^{\aleph_0}$. \qed
\end{fait}

Thus if $X$ is countable then $D^\alpha(X)=\emptyset$. In that
case we define the \emph{Cantor-Bendixson rank} of $x \in X$,
denoted $CB(x)$, by $CB(x)=\alpha$ if and only if $x \in
D^{\alpha}(X)\setminus D^{\alpha+1}(X)$.

Let us return  to our universal theory $T$. A \emph{complete
$n$-QF-type} (with respect to $T$), where the length $|\bar x|=n$,
is a set $p(\bar x)$ of quantifier-free formulas  such that for any
quantifier-free formula $\psi(\bar x)$ either $\psi \in p$ or
$\neg \psi \in p$, and such that there is a model $\cal M$ of $T$
having a tuple $\bar a$ such that $\mathcal M \models p(\bar a)$.

Let $S_n^{qf}(T)$ be the set of all completes $n$-QF-types of $T$. Then $S_n^{qf}(T)$ is equipped  with a topology as follows. Take
for basis open sets the sets of the form
$$
[\psi]=\{p \in S_n^{qf}(T)|\psi \in p\}, \hbox{ where $\psi$ is a
QF-formula}.
$$

The following fact is a consequence of the compactness theorem and the 
proof proceeds in a similar way to that of \cite[Th\'eor\`eme 4.04]{poizat}.

\begin{fait} \label{fact2}$S_n^{qf}(T)$ is a compact totally disconnected
space. \qed
\end{fait}

An \textit{$n$-marked} model of $T$ is a tuple  $(\mathcal M, \bar a)$, where
$\bar a$ is an $n$-generating tuple  of $\cal M$ and $\mathcal{M} \models T$.  For a  marked
model $(\mathcal M, \bar a)$ of $T$, we  let
$$
\Delta(\mathcal M, \bar a)=\{\psi(\bar x )\mid \mathcal M \models \psi(\bar a), \psi
\hbox{ is quantifier-free}\}. 
$$

It follows that   $\Delta(\mathcal M, \bar a) \in S_n^{qf}(T)$. 

Conversely,  let  $p(\bar x)
\in S_n^{qf}(T)$. By definition, there is  a model  $\mathcal N \models T$ having  a tuple $\bar
a$ such that $\mathcal N \models p(\bar a)$.  Then 
$(\mathcal M=\<\bar a\>, \bar a)$ is an $n$-marked model of
$T$ and $\Delta(\mathcal M, \bar a)=p(\bar x)$.  We notice that one can have  $\Delta(\mathcal M, \bar a)=\Delta(\mathcal M, \bar b)$, with $\bar a$  different from $\bar b$. 

Suppose now that $T$ has at most countably many finitely generated
models. Then $S_n^{qf}(T)$ is countable for any $n$, and by the
 Fact \ref{fact1}-\ref{fact2} one can associate an ordinal rank to marked models of
$T$ as follows.  Let $(\mathcal M, \bar a)$ be a marked model of $T$ where
$|\bar a|=n$. Then $\Delta(\mathcal M, \bar a) \in S_n^{qf}(T)$
and we define $Rk_n(\mathcal M, \bar a)=CB(\Delta(\mathcal M, \bar a))+1$, 
where the Cantor-Bendixson rank is calculated relatively to the topology defined above.

\begin{fait} \emph{\cite[Section 5]{ould-f-g}} Suppose that   $T$ has at most countably many finitely
generated models. If $\mathcal M$ is a finitely generated
model of $T$ and $\bar a$, $\bar b$ are finite generating tuples of
$\cal M$ with $|\bar a|=n$ and $|\bar b| =m$, then $Rk_n(\mathcal
M, \bar a)= Rk_m(\mathcal M, \bar b)$. 

\qed
\end{fait}

The above fact allows the definition of $Rk(\mathcal M)$ without
any reference to a particular generating set of $\cal M$.

If $G$ is an equationally noetherian group, then finitely generated $G$-limit groups coincide with finitely generated models of the universal theory of $G$ \cite[Theorem 2.1(4)]{ould-e-n}, and  there are  at most countably many of them 
\cite[Proposition 6.1]{ould-e-n}. Since a linear group $G$ is
equationally noetherian \cite[Theorem B1]{rem-mya-alge-geo},  we obtain that $Rk(L)$ is an ordinal for
any finitely generated $G$-limit group $L$. In fact, we have a bit more.

\begin{theo}\label{princip-theo} If $G$ is a linear group and $L$ is a finitely
generated model of the universal theory of $G$ \emph{(}i.e., $L$
is a finitely generated $G$-limit group\emph{)}, then $Rk(L)$ is finite.
\end{theo}

Since free groups are linear, it follows that the Cantor-Bendixson rank of a limit group is finite. Before proving the above theorem, we  explain the relation between the space of $n$-marked limit groups of
free groups $\mathcal F_n$, and  the space $S_n^{qf}(T)$ where
$T$ is the universal theory of $F_2$.  

Let us recall
the definition of $\mathcal F_n$  from
\cite{Cham-Guir}. Let $2^{ F_n}$ be the set of all subsets of the
free group $F_n$. Let us denote by $B_r$ the ball of radius $r$ in
the free group $F_n$. For subsets $X,Y \in 2^{ F_n}$, define the
valuation
$$
v(X,Y)= \max \{r \in \mathbb N \cup \{+ \infty\}| X \cap B_r=Y
\cap B_r\}.
$$

The space $\mathcal G_n$ is viewed as the set of normal subgroups
of $F_n$ with the metric $d(X,Y)=e^{-v(X,Y)}$ and $\mathcal F_n$
as the subspace of $\mathcal G_n$ of elements $X$ such that
$F_n/X$ is a limit group. Then $\mathcal G_n, \mathcal F_n$ are 
compact and totally disconnected.

For any $p \in S_n^{qf}(T)$, let $p^+$ be the set of
reduced words $w(\bar x)$ such that $w(\bar x)=1 \in p$.   Then
$p^+$ is a normal subgroup of $F_n$. 

\begin{fait} The function $f :  S_n^{qf}(T) \rightarrow \mathcal F_n$ defined by $f(p)=p^+$ is a homeomorphism. 
\end{fait}

\proof

Let us first prove that $f$ is a bijection.  Clearly $f$ is injective, because any element of  $S_n^{qf}(T)$  is a complete QF-type. Now let $X \in \mathcal F_n$. Then $H=F_n/X$ is a limit group, and thus $H$ is a model of the universal theory of $F_2$. We write $H=\<\bar a| w(\bar a)=1, w \in X\>$ and we let $p(\bar x)=\Delta(H, \bar a)$.  Clearly we have $f(p)=X$.

We are going to show that $f$ is countinous, the continuity of $f^{-1}$ can be treated in a similar way.  To this end it is sufficient to show that  for any $m \in \mathbb N$, for any  $X \in \mathcal F_n$, there is a QF-formula $\psi (\bar x)$ such that 
$$
f^{-1}(\{Y \in \mathcal F_n|v(X,Y)> m \})=[\psi]. 
$$

Let $P^+$ (resp. $P^{-}$)  be the set of reduced  words $w \in B_{m+1}$ such that  $H\models w(\bar a)=1$ (resp. $H \models w(\bar a) \neq 1)$, where  $H=\<\bar a| w(\bar a)=1, w \in X\>$. We let 
$$
\psi(\bar x) = \bigwedge_{w \in P^+} w(\bar x)=1\wedge \bigwedge_{v \in P^{-}} v(\bar x) \neq 1. 
$$

Now it is easily seen that $
f^{-1}(\{Y \in \mathcal F_n|v(X,Y)> m \})=[\psi]
$.  \qed

\bigskip

To prove Theorem \ref{princip-theo} we will use the fact that the Morley
rank of $GL_n(K)$, where $K$ is a pure  algebraically closed field, is
finite. We recall  some definitions.  The Cantor rank is defined for constructible sets in topological spaces.  The next definition is  more  suited to our context.

\begin{defi} \emph{(Cantor rank)}  Let $X$ be a set and let $\cal B$ be a boolean algebra in $X$. The Cantor rank
 $\rm RC(Y)$ of a  subset $Y \in \cal B$  is defined
by induction  as follows:

$\bullet$ $\Rc(Y) \geq 0$ if and only if $Y \neq \emptyset$;

$\bullet$ $\Rc(Y) \geq \alpha +1$ if and only if for every $n \in
\mathbb N$,  there exist a sequence of disjoints subsets $(Y_i \in
\mathcal{B}| i \leq n)$ of $X$ such that $\Rc(Y_i \cap Y) \geq
\alpha$;

$\bullet$ if $\alpha$ is a limit ordinal,  then $\Rc(Y) \geq \alpha$
if and only if $\Rc(Y) \geq \beta$ for every $\beta < \alpha$.

We put $\Rc(Y)$=\emph{sup} $\{\alpha ~|~\Rc(Y) \geq \alpha\}$
if $Y \neq \emptyset$; $\Rc(\emptyset)=-1$. \qed

\end{defi}

\begin{defi} \emph{(Morley rank for groups)} Let $G$ be a group and let $^*G$ be an
$\aleph_0$-saturated elementary extension of $G$.  Let $\mathcal B_n$ be the boolean algebra of
definable sets in $^*G^n$. We define the Morley rank of $X \in
\mathcal B_n$,  by $RM(X)=RC(X)$. We say that $G$ has a finite
Morley rank if $RM(\bar x=\bar x)$, $|\bar x|=n$,  is finite for
any $n$. We let $RM(G)=RM(x=x)$. \qed
\end{defi}

Now we introduce the QF-rank which is computed relatively to the boolean algebra of QF-definables sets. 

\begin{defi}\emph{(Quantifier-free rank for groups)}  Let $G$ be a group and let $^*G$ be an
$\aleph_0$-saturated elementary extension of $G$.  Let $\mathcal B_n$ be the boolean
algebra of quantifier-free definable sets in $^*G^n$. We define
the quantifier-free  rank of $X \in \mathcal B_n$ by
$QF(X)=RC(X)$. We say that $G$ has a finite QF-rank if $QF(\bar
x=\bar x)$, $|\bar x|=n$,  is finite for any $n$.    We
let $QF(G)=QF(x=x)$. \qed
\end{defi}

These definitions does not depend on a particular choice of
$^*G$.  The proof of that fact in the case of QF-definables sets, follows the same method as that of \cite[Lemma 6.2.2, Lemma 6.2.3]{marker}. We notice also that $^*G$ can be taken to be a nonprincipal ultrapower of $G$ (\cite[Theorem
6.1.1]{Chang-Keisler} or \cite[Exercice 4.5.37]{marker}).

\smallskip
If $\psi(\bar x)$ is a formula with parameters from $G$ (i.e., an
$\mathcal L_G$-formula), we define $RM(\psi)$ to be the Morley
rank of the set $\psi(G)=\{\bar a \in {^*G}|^*G \models \psi(\bar
a)\}$. If $\psi$ is QF-formula, $QF(\psi)$ is defined similarly.

\begin{fait}\emph{\cite[Corollary 6.2.23]{marker}}\label{fact-mor-algebraic} If $K$ is a pure  algebraically closed field then
$GL_n(K)$ has a finite Morley rank. Furthermore
$RM(GL_n(K))=\dim(GL_n(K))$, where $\dim$ denotes the dimension in
the usual Zariski topology. \qed
\end{fait}

In the following lemma, we write $RM_{\mathcal M}$,  $QF_{\mathcal M}$,  to distinguish the model $\mathcal M$ relatively to   which the rank is calculated.

\begin{lem} \label{lem-pricip}Let $\mathcal M$ and $\mathcal N$ be two groups such that $\mathcal M \leq \mathcal N$. Then for every quantifier-free
$\mathcal L_{\mathcal M}$-formula $\phi(\bar x)$, 
$QF_{\mathcal M}(\phi) \leq QF_{\mathcal N}(\phi)$. In particular
$QF(\mathcal{M}) \leq QF(\mathcal{N})$. We also have  $QF_{\mathcal
N}(\phi) \leq RM_{\mathcal N}(\phi)$.
\end{lem}

\proof

We claim first that there
exists an $\aleph_0$-saturated elementary extension
$\mathcal{M}_1$ of $\cal M$ and an $\aleph_0$-saturated elementary
extension $\mathcal{N}_1$ of $\cal N$ such that $\mathcal{M}_1
\leq \mathcal{N}_1$ and the following diagram commute
$$
\begin{array}{ccc}
  \mathcal{M}& \leq & \mathcal{N} \\
  \downarrow& \; & \downarrow \\
  \mathcal{M}_1 & \leq & \mathcal{N}_1 \\
\end{array}
$$
Let $\mathcal{M}_1$ be any $\aleph_0$-saturated elementary
extension of $\cal M$. Let 
$$
Diag(\mathcal{M}_1)=\{\phi| \phi \hbox{ is a quantifier-free }\mathcal L_{\mathcal M_1} \hbox{-formula and } \mathcal M_1 \models \phi\}, 
$$
$$
Diag_{el}(\mathcal{N})=\{\phi| \phi \hbox{ is an }\mathcal L_{\mathcal N} \hbox{-formula and } \mathcal N \models \phi\},
$$
and let $\Gamma=Diag(\mathcal{M}_1) \cup
Diag_{el}(\mathcal{N})$ in the language $\mathcal L_{\mathcal M_1}  \cup \mathcal L_{\mathcal N}$. 

Suppose that $\Gamma$ is not consistent.  Then there exists a
quantifier-free $\mathcal{L}_{\cal M}$-formula $\psi(\bar x)$ such
that $\mathcal{M}_1 \models \exists \bar x \psi(\bar x)$ and
Diag$_{el}(\mathcal{N})\models \forall \bar x \neg \psi(\bar x)$;
which is clearly a contradiction. Therefore $\Gamma$  is consistent and  we take
$\mathcal{N}_1$ to be any $\aleph_0$-saturated elementary
extension of a model of $\Gamma$. This ends the proof of our claim.

We  now prove by induction on $\alpha$ that for any
quantifier-free $\mathcal{L}_{\mathcal{M}_1}$-formula $\phi$,
$RC_{{\cal M}_1}(\phi) \geq \alpha \Rightarrow RC_{{\cal
N}_1}(\phi) \geq \alpha$, where here for a model $\mathcal K$, $RC_{{\cal K}}(\phi)$ denotes the Cantor rank computed in $\cal K$. 

Since $\phi$ is quantifier-free the result is clear for
$\alpha=0$. If $\alpha$ is a limit ordinal the result follows by 
induction.

For $\alpha+1$. Let $n \geq 1$. There are quantifier-free  $\mathcal{L}_{\mathcal
M_1}$-formulas $\psi_1(\bar x)$, $\psi_2(\bar x)$, $\dots$, $\psi_n(\bar x)$,  such
that $\psi_1(\mathcal {M}_1)$, $\psi_2(\mathcal {M}_1)$, $\dots$, $\psi_n(\mathcal M_1)$ is a family of pairwise disjoint subsets of
$\phi(\mathcal {M}_1)$ and $QF_{\mathcal M_1}(\psi_i(\bar x)) \geq
\alpha$ for all $i$. Define $\varphi_i(\bar x)$ as follows
$$
\varphi_i(\bar x)= \phi(\bar x) \wedge \psi_i(\bar x) \wedge_{j<i}
\neg \psi_j(\bar x). 
$$

Then $\mathcal{M}_1 \models \forall \bar x (\varphi_i(\bar x)
\Leftrightarrow \psi_i(\bar x))$. Therefore $RC_{\mathcal
M_1}(\psi_i(\bar x))=RC_{\mathcal M_1}(\varphi_i(\bar x))$ and by
induction we have $RC_{\mathcal N_1}(\varphi_i(\bar x)) \geq
\alpha$. Now $\varphi_1(\mathcal {N}_1)$, $\varphi_2(\mathcal
{N}_1)$, $\dots$, $\varphi_1(\mathcal {N}_1)$ is a  family of pairwise disjoint
subsets of $\phi(\mathcal {N}_1)$ and $RC_{\mathcal
N_1}(\varphi_i(\bar x)) \geq \alpha$ for all $i$. 

Thus
$RC_{\mathcal N_1}(\phi(\bar x)) \geq \alpha+1$.

The proof of $QF_{\mathcal N}(\phi) \leq RM_{\mathcal N}(\phi)$ proceeds in a similar way.  \qed

Combining Fact \ref{fact-mor-algebraic} and Lemma \ref{lem-pricip}, we get the following. 

\begin{cor} A linear group has a finite QF-rank. \qed
\end{cor}

By the precedent corollary, Theorem \ref{princip-theo} is  a consequence of the
following more general proposition which is analogous to the one in the full category of definable sets \cite[Exercice 6.6.19]{marker}. 

\begin{prop} Let $G$ be a group of finite QF-rank. Then for any
finitely generated model $L$ of the universal theory of $G$,
$Rk(L)$ is finite.
\end{prop}

\proof

Let $L$ be a finitely generated model of the universal theory $T$
of $G$, $n$-generated by the tuple $\bar a$. Since $Rk(L)=Rk_n(L,
\bar a)=CB(\Delta(L, \bar a))+1$, it is sufficient to show that for
any $p \in S_n^{qf}(T)$, $CB(p)$ is finite. 

Let $X=S_n^{qf}(T)$. Let $ QF_0(p)= \hbox{ min } \{QF(\psi)\mid \psi \in p \}$.  Then
$QF_0(p)$ is a natural number as $G$ has a finite QF-rank. We are
going to prove by induction on natural numbers $m \geq 0$, that if
$CB(p)=m$, then $QF_0(p) \geq m$.

For $m=0$. Since $CB(p)=0$, $p \in X\setminus  D(X)$ and thus $p$
is isolated. Therefore,  there exists a QF-formula $\psi_0(x_1,
\cdots, x_n )$ such that
$$
T \vdash \forall \bar x(\psi_0(\bar x )\Rightarrow \psi(\bar x)),
$$
for any $\psi \in p$. Therefore $ \psi_0 \in p$ and $QF_0(p)=QF(\psi_0)$. Since
$\exists \bar x \psi_0(\bar x)$ is consistent with $T$,  we get
$QF(\psi_0) \geq 0$ as desired.

For $m+1$. Then $p \in D^{m+1}(X)\setminus  D^{m+2}(X)$. Therefore
there exists a QF-formula $\psi(\bar x)$ such that $\psi \in p$
and for any $q \in D^{m+1}(X)$ if $\psi \in q$ then $p=q$. We may
assume that $QF_0(p)=QF(\psi)$.

Since $p$ is not isolated in $D^m(X)$, there exists $q \in D^m(X)$
such that $\psi \in q_1$ and clearly we must have $q_1 \in
D^{m}(X)\setminus D^{m+1}(X)$. Since $q_1 \neq p$ there is a
QF-formula $\psi_1 \in p$ such that $\neg \psi_1 \in q_1$. By the
same argument, there exists $q_2  \in D^{m}(X)\setminus
D^{m+1}(X)$ and $\psi_2$ such that $\psi \wedge \psi_1 \in q_2$,
$\psi_2 \in p$, $\neg \psi_2 \in q_2$. Containing in this manner,
we get an infinite sequence $((q_n, \psi_n)|n \in \mathbb N)$,
such that $q_n \in D^{m}(X)\setminus D^{m+1}(X)$, $\psi_n \in p$,
$\neg \psi_n \in q_n$, $\psi_1, \cdots, \psi_{n-1} \in q_n$.

By induction $QF(\psi_n) \geq m$, as $CB(q_{n+1}) =m$ and
$\psi_{n} \in q_{n+1}$. We also have $QF(\psi \wedge \psi_n) \geq
m$ and $QF(\psi \wedge \neg \psi_n) \geq m$. Then the QF-formulas
$$
\psi \wedge \neg \psi_1, \psi \wedge \psi_1 \wedge \neg \psi_2,
\psi \wedge \psi_1 \wedge \psi_2, \psi \wedge \psi_1 \wedge \psi_2
\wedge \neg \psi_3, \cdots
$$

are disjoint QF-formulas each of which has a QF-rank greater than
$m$. Therefore $QF(\psi) \geq m+1$. \qed

We end this note with the following lemma.

\begin{lem}$\;$

$(1)$ Let $T$ be the universal theory of the free group $F_2$ of
rank 2 and let $L$ be a finitely generated model of $T$. Then
$Rk(L)=1$ if and only if $L\cong 1$ or $L \cong \mathbb Z$ or $L
\cong F_2$.

$(2)$ $Rk(\mathbb Z^n)=n$, $Rk(F_n) \geq n-1$ for any $n \geq 2$.
\end{lem}

\proof

$(1)$  If $Rk(L)=1$, then  $L$ is embeddable in $F_2$ \cite[Lemma 6.5]{ould-f-g} and
thus $L$ is free. It is also clear that $Rk(1)=Rk(\mathbb Z)=1$, because $1$ and $\mathbb Z$  are isolated by the formula $x=1$ and $x \neq 1$
respectively.  Since any subgroup of a limit group, generated by
two elements is either abelian or free of rank $2$ \cite[Lemma 5.5.4]{chiswell}, it follows that
$F_2$ is isolated by the formula $[x,y] \neq 1$.

Let $F_n$ be a free group of rank $\geq 3$ and let $\{x_1, \cdots,
x_n\}$ be a basis of $F_n$. Then $F_{n-1} \preceq_{\exists}
F_{n-1}*\<x_n|\>$(see \cite{sel1}). Therefore if $Rk(F_n)=1$ and if $\psi(x_1,
\cdots, x_n)$ is a QF-formula isolating $F_n$, then there exists
$b \in F_{n-1}$ such that $F_{n-1} \models \psi(x_1, \cdots,
x_{n-1},b)$ and thus we find $F_n \cong F_{n-1}$, which is a
contradiction. Therefore $Rk(F_n) \geq 2$.

$(2)$ Similar to the above proof. \qed

\bigskip
\noindent \textbf{Acknowledgements.} We thank B. Poizat and F. Wagner   for having pointed out that the QF-rank is bounded by the Morley rank.

\noindent Abderezak OULD HOUCINE,\\
Universit\'e de Lyon;\\
Universit\'e Lyon 1;\\
INSA de Lyon, F-69621;\\
Ecole Centrale de Lyon;\\
CNRS, UMR5208, Institut Camille Jordan,\\
43 blvd du 11 novembre 1918,\\
F-69622 Villeurbanne-Cedex, France\\
E-mail address : \ttfamily{ould@math.univ-lyon1.fr} 
\end{document}